\documentclass{IEEEtran4PSCC}

\ifCLASSINFOpdf
   \usepackage[pdftex]{graphicx}
  
\else
 
   \usepackage[dvips]{graphicx}
  
\fi
\usepackage[cmex10]{amsmath}
\ifCLASSOPTIONcompsoc
  \usepackage[caption=false,font=normalsize,labelfont=sf,textfont=sf]{subfig}
\else
  \usepackage[caption=false,font=footnotesize]{subfig}
\fi
\usepackage{url}
\usepackage{xcolor}
\usepackage{mathtools, amsfonts, amssymb}
\usepackage{hyperref}
\usepackage{caption}
\usepackage{subcaption}
\usepackage{subfig}

\newcommand{\by}{\boldsymbol{y}}
\newcommand{\bx}{\boldsymbol{x}}
\begin{document}
\title{Hierarchical Graph Modeling for\\Multi-Scale Optimization of Power Systems}

%% To specify the authors when (number of affiliations <= 2)
%\author{
%\IEEEauthorblockN{David L. Cole\\ Victor M. Zavala}
%\IEEEauthorblockA{Department of Chemical and Biological Engineering \\
%University of Wisconsin-Madison\\
%Madison, Wisconsin, United States\\
%\{dlcole3, zavalatejeda\}@wisc.edu}
%\and
%\IEEEauthorblockN{Harsha Gangammanavar}
%\IEEEauthorblockA{Operations Research and Engineering Management \\
%Southern Methodist University\\
%Dallas, Texas, United States\\
%harsha@smu.edu}
%}

%% To specify the authors when (number of affiliations > 2)
 \author{\IEEEauthorblockN{David L. Cole\IEEEauthorrefmark{1},
 Harsha Gangammanavar\IEEEauthorrefmark{2},
 Victor M. Zavala\IEEEauthorrefmark{1}\IEEEauthorrefmark{3}, 
}
 \IEEEauthorblockA{\IEEEauthorrefmark{1} Department of Chemical and Biological Engineering, 
 University of Wisconsin-Madison,
 Madison, WI, United States\\
 dlcole3@wisc.edu}
 \IEEEauthorblockA{\IEEEauthorrefmark{2} Operations Research and Engineering Management, 
 Southern Methodist University, Dallas, TX, United States\\
 harsha@smu.edu}
 \IEEEauthorblockA{\IEEEauthorrefmark{3} Mathematics and Computer Science Division, 
 Argonne National Laboratory, Lemont, IL, United States\\
 zavalatejeda@wisc.edu}
 }

\maketitle

\begin{abstract}
Hierarchical optimization architectures are used in power systems to manage disturbances and phenomena that arise at multiple spatial and temporal scales. We present a graph modeling abstraction for representing such architectures and an implementation in the {\tt Julia} package {\tt Plasmo.jl}. We apply this framework to a tri-level hierarchical framework arising in wholesale market operations that involves day-ahead unit commitment, short-term unit commitment, and economic dispatch. We show that graph abstractions facilitate the construction, visualization, and solution of these complex problems. 
\end{abstract}


\vspace{-0.2in}
\section{Introduction}
Hierarchical optimization architectures are used in power systems (and many other industrial systems) for managing operations, disturbances, and phenomena that arise at multiple spatial and temporal scales. These architectures involve multiple decision-making layers where decisions of higher layers influence or inform lower layers (and vice versa); for example, market operations often involve the solution of a unit commitment (UC) problem whose solution informs an economic dispatch (ED) problem \cite{conejo2018}.  Hierarchical decomposition is often necessary for enabling scalable implementation (e.g., solving a combined UC/ED problem in real-time might be impossible) and for providing intuitive decomposition of functionalities (which can aid explainability).  Capturing the unique characteristics of optimization problems in different hierarchical layers (i.e., space/time resolution, data, variables, objectives, constraints) and their hierarchical coupling is essential for enabling decision-making consistency across scales. This has motivated research in models and solution approaches that aim to identify how to best design hierarchical architectures to manage diverse types of features (e.g., identify the number of layers, resolutions, and decisions made by each layer). For example, Atakan and co-workers \cite{atakan2022} presented a stochastic optimization framework that consists of a tri-level hierarchy of market operations (day-ahead UC, short-term UC, and ED) that aims to handle high renewable penetration. The authors demonstrate that the hierarchical framework provides significant operational improvement over competing architectures. Guo and co-workers \cite{guo2017hierarchical} used a hierarchical architecture for decentralizing ED of a large power networks; this was a tri-level architecture where local, clustered agents (lowest layer) inform leader agents (middle layer), which in turn inform a coordinating agent (top layer).  Kong and co-workers \cite{kong2017hierarchical} proposed a hierarchical architecture for a network of electric vehicle charging stations connected to the grid; the formulation considers the placement of stations, the allocation of resources, and the operation policy of the stations on three separate hierarchical layers. They found that the framework provided improved system performance and quality of service. 

As power systems become increasingly complex (e.g., they include new assets and face new disturbances), it will become necessary to have modeling and solution tools that enable the seamless construction, evaluation, and benchmarking of different hierarchical architectures. In this work, we propose a graph-based modeling framework for representing hierarchical optimization structures arising in power system operations. The use of graphs to model structured optimization problems has been recently explored \cite{jalving2022graph, cole2022julia, berger2021gboml_tutorial, berger2021remote, allman2019}. A variety of tools for exploiting graph and graph-like structures are also available in open-source packages such as Plasmo.jl (in Julia) \cite{jalving2022graph} and Pyomo (in Python)  \cite{hart2011pyomo,hart2017pyomo}. In this work, we focus on the use of {\tt Plasmo.jl}; this package uses an {\tt OptiGraph} abstraction, where nodes of the graph contain optimization subproblems (with their own objective functions, data, variables, and constraints) and where edges capture connectivity (constraints) across subproblems. The {\tt OptiGraph} abstraction is flexible in that nodes can contain subproblems of different granularity; moreover,  the abstraction enables the creation of hierarchical structures (a node can be a graph itself). The graph abstraction provides the ability to build complex structures in a modular manner (e.g., node by node) and the ability to visualize, decompose, and aggregate the overall problem graph. We provide a case study show how the graph representation can be used for expressing and solving complex hierarchical problems arising in power systems.

\section{Graph-Based Modeling Overview}

\begin{figure}[!htp]
    \centering
    \vspace{-0.25in}
    \includegraphics[width=2.5in]{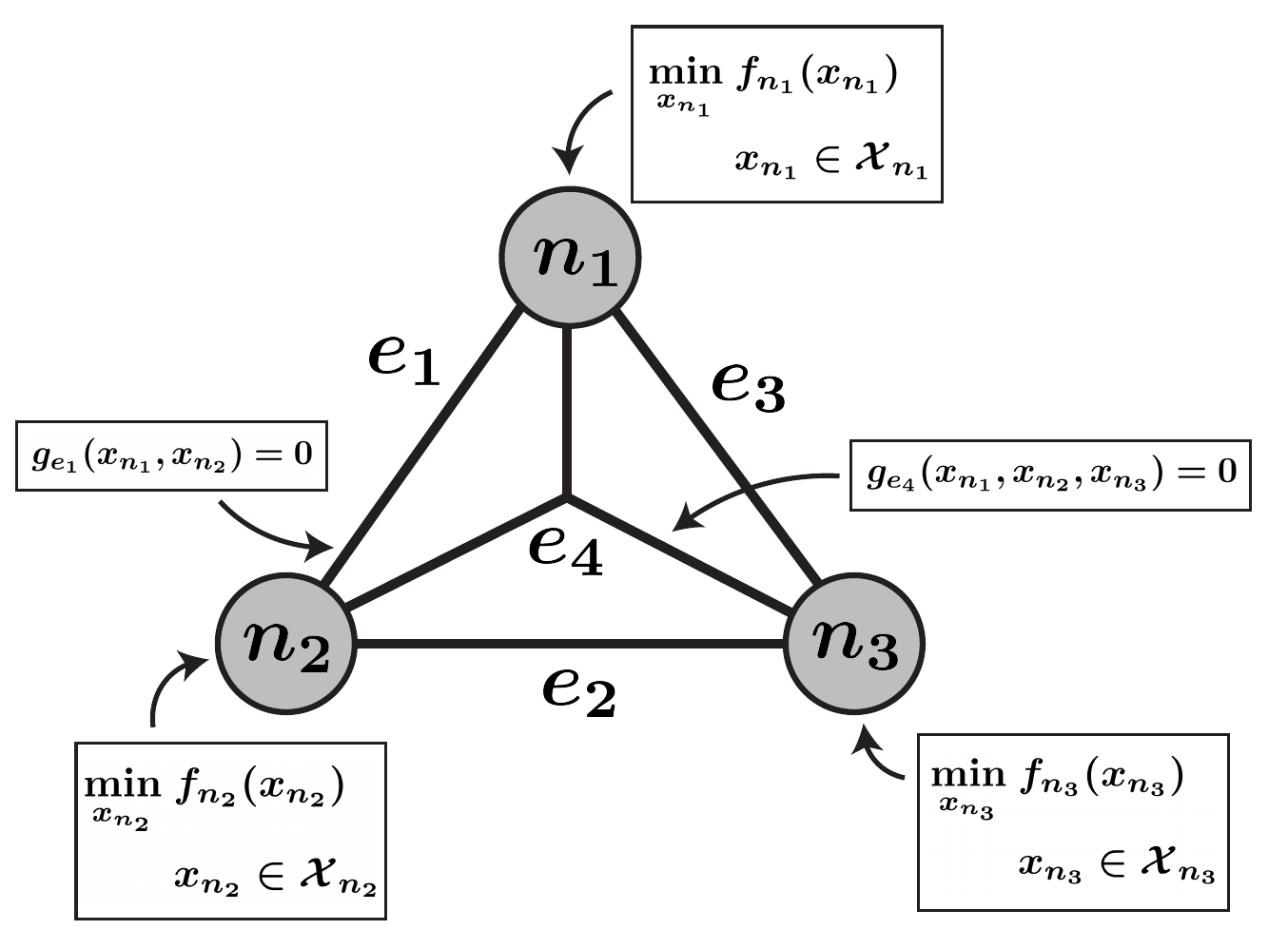}
    \caption{{\tt OptiGraph} abstraction used in {\tt Plasmo.jl}.}
    \label{fig:plasmo_overview}
    \end{figure}

{\tt Plasmo.jl} is a Julia package that models general optimization problems as hypergraphs. This package has been described in detail by Jalving and co-workers \cite{jalving2022graph}, but here we provide a short overview of how this can be used for representing hierarchical problems. {\tt Plasmo.jl} is built on an abstraction called {\tt OptiGraphs}, which are graphs containing {\tt OptiNodes} ($\mathcal{N}$) and {\tt OptiEdges} ($\mathcal{E}$). {\tt OptiNodes} contain subproblems (with their own variables, constraints, data, and objective functions), and {\tt OptiEdges} are linking constraints that capture connectivity between {\tt Optinodes}. We denote an {\tt OptiGraph} as $\mathcal{G}(\mathcal{N}, \mathcal{E})$, where $\mathcal{N}(\mathcal{G})$ is the set of {\tt OptiNodes} in $\mathcal{G}$ and $\mathcal{E}(\mathcal{G})$ is the set of {\tt OptiEdges} in $\mathcal{G}$. A visualization of an {\tt OptiGraph} containing three {\tt OptiNodes} is shown in Figure \ref{fig:plasmo_overview}. The optimization model associated with an {\tt OptiGraph} can be represented as
\begin{align}\label{eq:optigraph}
    \begin{split}
    \min_{\{ x_n \}_{n \in \mathcal{N}(\mathcal{G})}} &\; \sum_{n \in \mathcal{N}(\mathcal{G})} f_n(x_n)\\
    \textrm{s.t.} &\; x_n \in \mathcal{X}_n, \quad n \in \mathcal{N}(\mathcal{G})\\
    &\; g_e(\{x_n\}_{n \in \mathcal{N}(e)}) \geq 0, \quad e \in \mathcal{E}(\mathcal{G})
    \end{split}
\end{align}
where $\mathcal{N}(e)$ is the set of {\tt OptiNodes} that support {\tt OptiEdge} $e$. The notion of nodes and edges is highly flexible in this abstraction; for instance, in a power system context, a node can represent a spatial location, time instance, a specific asset, or an entire network.  Moreover, each node can have its own independent features (e.g., data, objective functions, constraints). The edges (containing the constraints $g_e$) can be used to link nodes across time, space, or hierarchical layers.

{\tt OptiGraphs} enable hierarchical representations and modular model building via the use of subgraphs. Specifically, within {\tt Plasmo.jl}, an {\tt OptiGraph} can be embedded in another {\tt OptiGraph} as a node. For example, consider the {\tt OptiGraphs} $\mathcal{G}_i$ and $\mathcal{G}_j$, each with an independent set of {\tt OptiNodes} and {\tt OptiEdges}. We can consider these {\tt OptiGraphs} as low-level graphs (also referred to as subgraphs) that can be used to build a higher-level {\tt OptiGraph} which we denote by $\mathcal{G}(\{ \mathcal{G}_i, \mathcal{G}_j \}, \mathcal{N}_g, \mathcal{E}_g)$. The set of nodes $\mathcal{N}_g$ are contained on $\mathcal{G}$ and are separate from $\mathcal{N}(\mathcal{G}_i)$ and $\mathcal{N}(\mathcal{G}_j)$, such that $\mathcal{N}_g = \mathcal{N}(\mathcal{G}) / \{\mathcal{N}(\mathcal{G}_i) \cup \mathcal{N}(\mathcal{G}_j)\}$. Similarly, $\mathcal{E}_g = \mathcal{E}(\mathcal{G}) / \{ \mathcal{E}(\mathcal{G}_i) \cup \mathcal{E}(\mathcal{G}_j) \}$, meaning $\mathcal{E}_g$ may connect nodes across $\mathcal{N}_g$, $\mathcal{N}(\mathcal{G}_i)$, and/or $\mathcal{N}(\mathcal{G}_j$). The {\tt OptiGraph} $\mathcal{G}$ may also be placed in another higher-level {\tt OptiGraph} $\mathcal{G}'(\{\mathcal{G}\}, \mathcal{N}_g', \mathcal{E}_g' )$. Any subgraph can also be collapsed into a single {\tt OptiNode} (containing the entire problem of the subgraph); this is useful for visualizing hierarchies. For instance, the hierarchical setting allows us to capture how assets can be aggregated in space (e.g., assets can be embedded at a network location), or how multiple time points can be embedded in another time point (e.g., multiple 5-min time periods can be embedded in an hour). This feature is key for representing hierarchical structures that span multiple scales.

    \begin{figure}[!htp]
\centering
\vspace{-0.2in}
\includegraphics[width=1.5in]{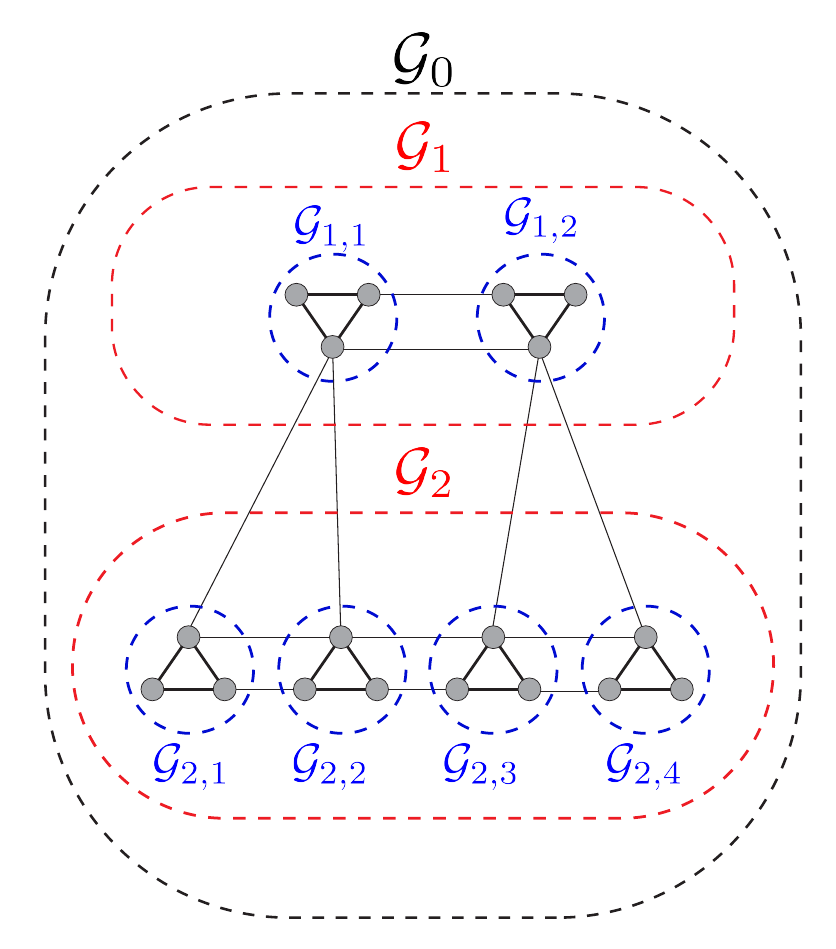}
\caption{{\tt OptiGraph} abstraction of bi-level problem.}
\label{fig:2_layer_prob}
\end{figure}
The general approach for representing hierarchical problems as graphs is illustrated in Figure \ref{fig:2_layer_prob}. This is a bi-level hierarchical problem; the top-level {\tt OptiGraph} is given by $\mathcal{G}_1(\{ \mathcal{G}_{1,1}, \mathcal{G}_{1, 2}\}, \emptyset, \mathcal{E}_1)$, the lower-level {\tt OptiGraph} by $\mathcal{G}_2(\{ \mathcal{G}_{2, 1}, \mathcal{G}_{2,2}, \mathcal{G}_{2, 3}, \mathcal{G}_{2, 4}\}, \emptyset, \mathcal{E}_2)$, and the overall {\tt OptiGraph} by $\mathcal{G}_0(\{\mathcal{G}_1, \mathcal{G}_2\}, \emptyset, \mathcal{E}_0)$, where $\mathcal{E}_0$ are the constraints linking the solutions of the upper and lower layers. 

{\tt OptiGraphs} enable flexible partitioning of hierarchical structures, and this can be used to implement different solution approaches. For instance, it has been recently shown that graph structures facilitate the development of decomposition algorithms \cite{jalving2022graph,cole2022julia, allman2019,daoutidis2019decomposition, shin2021graph, shin2020decentralized}. To be specific, any {\tt OptiNode} or subgraph can be treated as an individual optimization problem; for example, the {\tt OptiGraph} presented in Figure \ref{fig:2_layer_prob} can be solved in at least three different ways (each likely resulting in different solutions): i) $\mathcal{G}_0$ could be solved as a single monolithic problem; ii) $\mathcal{G}_1$ and $\mathcal{G}_2$ can be solved sequentially with the solution of $\mathcal{G}_1$ passed via $\mathcal{E}_0$; iii) $\mathcal{G}_{1, 1}$ and $\mathcal{G}_{1, 2}$ can be solved sequentially with the solution of $\mathcal{G}_{1, 1}$ passed via $\mathcal{E}_1$. Then those solutions to $\mathcal{G}_1$ can be passed via $\mathcal{E}_0$ to $\mathcal{G}_2$, and $\mathcal{G}_{2, 1}$, $\mathcal{G}_{2, 2}$, $\mathcal{G}_{2, 3}$, and $\mathcal{G}_{2, 4}$ can be solved sequentially, with solutions passed via $\mathcal{E}_2$. We can thus see that the {\tt OptiGraph} abstraction offers significant flexibility in modeling and solving hierarchical problems. 

\section{Case Study}

\subsection{Problem Overview}

\begin{figure}[!ht]
    \centering
    \vspace{-0.1in}
    \includegraphics[width=3.0in]{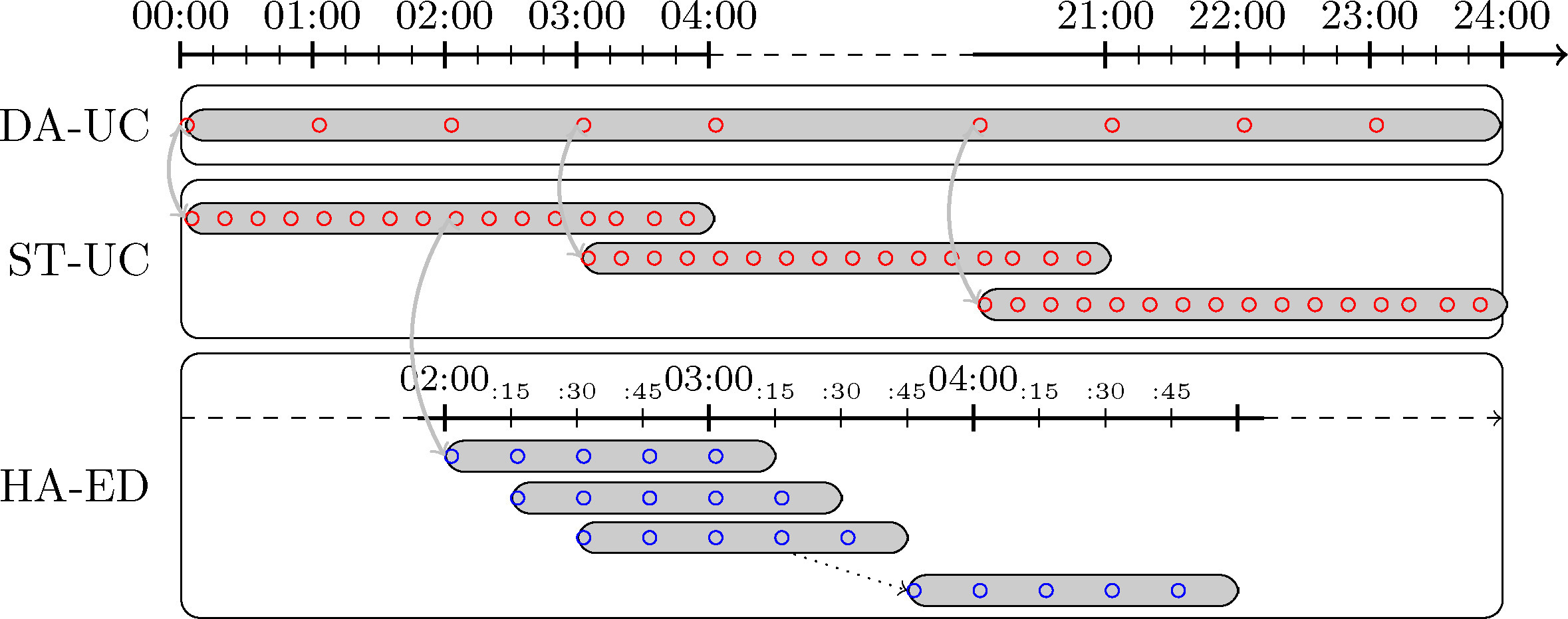}
    \caption{The tri-level hierarchical architecture of Atakan et al. \cite{atakan2022} (Reproduced with permission from Elsevier).} 
    \label{fig:HP_framework}    
\end{figure}

We consider the tri-level problem proposed in \cite{atakan2022} for capturing coupling in market operations (see Figure \ref{fig:HP_framework}). Each layer is composed of subproblems at different timescales and these are linked to subproblems in other layers. The top layer is a day-ahead unit commitment (DA-UC) problem that schedules a subset of conventional (non-renewable) generators (denoted as $\Gamma^d_c$). The DA-UC layer has a 1-hour resolution and a 24-hour horizon, and an entire time horizon is partitioned into periods of 24 hours (DA-UC is solved every 24 hours). The second layer includes a short-term unit commitment (ST-UC) problem; this schedules a subset of conventional generators (denoted as $\Gamma_c^s$), such that $\Gamma_c^d \cap \Gamma_c^s = \emptyset$, while also incorporating the commitment decisions of the DA-UC subproblems. The ST-UC layer has a 15-min resolution with subproblems containing a 4-hour horizons and solved every 3 hours (there is overlap). The bottom  layer is an hour-ahead economic dispatch (HA-ED) layer which determines the generation levels for units committed in the DA-UC and ST-UC layers. The HA-ED subproblems have a 15-minute resolution and a 75-min time horizon, and are solved every 15 minutes (there is overlap). Thus, for a given day, there are: 1 DA-UC subproblem, 8 ST-UC subproblems, and 96 HA-ED subproblems (12 for each ST-UC subproblem). We highlight that this architecture is just one design (of many possible ones). In other words, one could design diverse hierarchical architectures (e.g., experimenting with the types of variables, resolutions, and time horizons that each layer uses). 

The detailed model can be found in \cite{atakan2022}; here, we provide a high-level perspective to illustrate how complex models are embedded in the different layers and how coupling arises between layers. We use the sets $\Gamma$ for the set of all generators, $\Gamma_r$ for the set of renewable generators, and define $\Gamma_c^h = \Gamma_c^d \cup \Gamma_c^s$. We also use $*$ and $**$ to define variables, sets, or functions corresponding to a layer; here, the symbols $*$ or $**$ are exchanged for $d$, $s$, or $h$ to denote the DA-UC, ST-UC, or HA-ED layers, respectively. As each subproblem considers different sets of times, we use $\mathcal{T}^*_i$ for the set of times (in hours) of the $i$th subproblem of the $*$ layer. We also define the sets $\bar{\mathcal{T}}^*_i$ as the set of times without the first time point of the subproblem. We define $\Delta^*$ as the time step for the subproblem in hours($\Delta^d = 1$, $\Delta^s = 0.25$, and $\Delta^h = 0.25$). Note that we set $\Delta^s = \Delta^h$, and this has a small influence on the formulation of the problem presented below. The decision variables are given by:
\begin{align*}
    \begin{split}
    \bx^d_{i}  =&\; \big(x_{g, t}, s_{g, t}, z_{g, t}\big)_{\forall g \in \Gamma^d_c, t \in \mathcal{T}^d_i} \\
    \by^d_{i} =&\;  \Big((G^{+, d}_{g, t}, G^{-, d}_{g, t})_{\forall{g \in \Gamma^d_c \cup \Gamma_r}}, (F^d_{j, k, t})_{\forall (j, k)\in \mathcal{L}}, \\
    &\; (D^{d}_{j, t}, \theta^d_{j, t})_{\forall j \in \mathcal{B}}\Big)_{\forall t \in \mathcal{T}^d_i}\\
    \color{red} \bx^s_{i} =&\; \color{red}  \big(x_{g, t}, s_{g, t}, z_{g, t}\big)_{\forall g \in \Gamma^s_c, t \in \mathcal{T}^s_i} \\
    \color{red} \by^s_{i} =&\; \color{red} \Big((G^{+, s}_{g, t}, G^{-, s}_{g, t})_{\forall{g \in \Gamma_c^h \cup \Gamma_r}}, (F^s_{j, k, t})_{\forall (j, k)\in \mathcal{L}},\\
    &\; \color{red} \qquad (D^{s}_{j, t}, \theta^s_{j, t})_{\forall j \in \mathcal{B}}\Big)_{\forall t \in \mathcal{T}^s_i}\\
    \color{blue} \by^h_{i} =&\; \color{blue}  \Big((G^{+, h}_{g, t}, G^{-, h}_{g, t})_{\forall{g \in \Gamma_c^h \cup \Gamma_r}}, (F^h_{j, k, t})_{\forall (j, k)\in \mathcal{L}}, \\
    &\; \color{blue} \qquad (D^{h}_{j, t}, \theta^h_{j, t})_{\forall j \in \mathcal{B}}\Big)_{\forall t \in \mathcal{T}^s_i}
    \end{split}
\end{align*}
Symbols $\bx^d_i$ and $\by^d_i$ are binary and continuous decision variables for the $i$th subproblem of DA-UC, $\color{red} \bx^s_i$ and $\color{red} \by^s_i$ are decision variables for the $i$th subproblem of ST-UC, and $\color{blue} \by^h_i$ are decision variables for the $i$th subproblem o HA-ED.  Symbols $x_{g, t}$, $s_{g, t}$, and $z_{g,t}$ are binary variables indicating for time $t$ whether generator $g$ is on/off, was turned on, or was turned off. $G^{+, *}_{g, t}$ is the power of generator $g$ consumed by the grid at time $t$ and $G^{-, *}_{g, t}$ is the power overgenerated (for conventional generators) or curtailed (for renewable generators) from $g$ at time $t$. $F^*_{j, k, t}$ is the power flow of transmission line $(j, k)$ from bus $j$ to bus $k$ during time $t$. $D^*_{j, t}$ is the amount of load shed at bus $j$ for time $t$, and $\theta^*_{j, t}$ is the bus angle for bus $j$ at time $t$. We also use black, red, and blue color to denote variables for DA-UC, ST-UC, and HA-ED, respectively.

The objective functions in the different layers are comprised of a UC part, $f_u$, and an ED part, $f_e$:
\begin{align}\label{eq:f_u}
    f_{u}(\bx^*_i) := \sum_{t \in \mathcal{T}^*_i} \sum_{g \in \Gamma^*_c} \phi_g^s s_{g, t} + \phi_g^f x_{g, t} \Delta^*
\end{align}
\begin{align}\label{eq:f_e}
    \begin{split}
        f_{e}(\by^*_i) := &\; \sum_{t \in \mathcal{T}^*_i} \Big( \sum_{g \in \Gamma^*_c} \phi^v_g G^{+,*}_{g, t} + \sum_{j \in \mathcal{B}} \Big(\phi^u_j D^{*}_{j, t} + \\
        &\; \sum_{g \in \Gamma_j \cap \Gamma^*_c}\phi^o_g G^{-, *}_{g, t} + \sum_{g \in \Gamma_j \cap \Gamma_r} \phi^c_g G^{-, *}_{g, t}   \Big) \Big)
    \end{split}
\end{align}
The function $f_u$ accounts for the startup cost $\phi^s_g$ and the no-load cost $\phi^f_g$ for generator $g$ at time $t$. The no-load cost is multiplied by $\Delta^*$ since the DA-UC and ST-UC levels have different time resolutions. The function $f_e$ accounts for the variable cost, $\phi^v_g$, of the energy consumed by the grid, the cost of overgeneration, $\phi^o_g$, and the cost of curtailment, $\phi^c_g$, for generator $g$ at time $t$. It also accounts for the cost of unmet demand $\phi^u_j$. We next define constraints for the layers: 
\begin{align}\label{eq:c_d}
    \begin{split}
        c_d(\by_i^*) := &\; \Big( \sum_{j \in \mathcal{B}:(j, k) \in \mathcal{L}} F^*_{j, k, t} - \sum_{j \in \mathcal{B}:(k, j) \in \mathcal{L}} F^*_{k, j, t} + \\
        &\; \sum_{g \in \Gamma_k} G^{+, *}_{g, t} + D^{*}_{k, t} - \hat{D}^*_{k, t} - \hat{R}^*_{k, t} \Big)_{k \in \mathcal{B}, t \in \mathcal{T}_i^*}
    \end{split}
\end{align}
This requires that, at each time point, the flows coming into the bus, the power consumed by the grid, and the amount of unmet demand is equal to the demand, $\hat{D}^*_{k, t}$, and the reserve requirements, $\hat{R}^*_{k, t}$ for $k \in \mathcal{B}$. The power flow equations use a DC approximation: 
\begin{align}\label{eq:c_f}
    c_f(\by^*_{i}) := \left( F^*_{j, k, t} - B_{j, k} (\theta^*_{j, t} - \theta^*_{k, t}) \right)_{(j, k) \in \mathcal{L}, t \in \mathcal{T}^*_i}
\end{align}
The renewable resources are also restricted to a specific value; this is enforced by constraint $c_r$, where $\hat{G}^*_{g, t}$ is the amount of power produced by renewable generator $g$ at time $t$: 
\begin{align}\label{eq:c_r}
    c_r(\by^*_{i}) := \left( G^{+, *}_{g, t} + G^{-, *}_{g, t} - \hat{G}^*_{g, t} \right)_{g \in \Gamma_r, t \in \mathcal{T}^*_i}
\end{align}
The ramp-up and ramp-down constraints ($c_{ru}$ and $c_{rd}$) have different forms because DA-UC and ST-UC/HA-ED have different time resolutions. Ramping constraints containing start-up and shut-down constraints are:
\begin{align}\label{eq:c_su}
    \begin{split}
    c_{su}(\by^*_{i}, \bx^{**}_j) := &\;\Big( G^{+, *}_{g, t} + G^{-, *}_{g, t} - G^{+, *}_{g, t - \Delta^*} - G^{-, *}_{g, t - \Delta^*} -\\
    &\; (\overline{S}_g - \overline{R}_g \Delta^{**} - \underline{C}_g) s^{**}_{g, t + \Delta^{**}} - \\
    &\; (\overline{R}_g \Delta^{**} + \underline{C}_g) x^{**}_{g, t + \Delta^{**}} + \underline{C}_{g} x^{**}_{g, t} \Big)_{g \in \Gamma_c^{**}, t \in \Theta}
    \end{split}
\end{align}
\vspace{-0.3in}
\begin{align}\label{eq:c_sd}
    \begin{split}
    c_{sd}(\by^*_{i}, \bx^{**}_j) :=&\; \Big( G^{+, *}_{g, t - \Delta^*} + G^{-, *}_{g, t - \Delta^*} - G^{+, *}_{g, t} - G^{-, *}_{g, t} - \\
    &\; (\underline{S}_g - \underline{R}_g \Delta^* - \underline{C}_g) z^{**}_{g, t + \Delta^{**}} -\\
    &\; (\underline{R}_g \Delta^{**} + \underline{C}_g) x^{**}_{g, t} +  \underline{C}_{g} x^{**}_{g, t + \Delta^{**}} \Big)_{g \in \Gamma_c^{**}, t \in \Theta}
    \end{split}
\end{align}
Here, $\overline{S}_g$ and $\underline{S}_g$ are the startup/shutdown limits for $g$, $\overline{R}_g$ and $\underline{R}_g$ are the ramp-up and ramp-down limits for $g$ as a function of time, and $\underline{C}_g$ is the minimum capacity of $g$. The set $\Theta$ is defined in this context as $\Theta = \{ t_i: t_i - \Delta^{*}= t_j - \Delta^{**}, \forall t_i \in \bar{\mathcal{T}}^*_i, t_j \in \mathcal{T}_j^{**} \}$. The model uses constraints $c_{su, 0}(\by^*_i, \by^h_k, \bx_j^{**}, \bx_{j - 1}^{**})$ and $c_{sd, 0}(\by^*_i,\by^h_k, \bx_j^{**}, \bx_{j - 1}^{**})$ to link the first time point of the $i$th subproblem with the solutions of the previoussubproblems (this introduces complex time coupling).  We also define operating regions:
\begin{align}
    \mathcal{X}_i^* &:=  \left\{ \bx^*_i | (x_{g,t}, s_{g, t}, z_{g, t}) \in \{0, 1\}^3, \quad \substack{\forall g \in \Gamma_c^*,\\ t \in \mathcal{T}_i^*} \right\}\label{eq:x}\\    
    \mathcal{Y}^* &:= \left\{ \begin{array}{c|l} & \underline{\theta}_j \le \theta_{j, t}^d \le \overline{\theta}_j, \quad \forall j \in \mathcal{B}, t \in \mathcal{T}_i^d\\
        & \underline{F}_{j, k} \le F^*_{j, k, t} \le \overline{F}_{j, k}, \quad \forall (j, k) \in \mathcal{L}, t \in \mathcal{T}_i^d\\
        \by^*_i  & (G^{+, *}_{g, t}, G^{-, *}_{g, t}) \in \mathbb{R}_{+}^2, \quad \forall g \in \Gamma_c^d \cup \Gamma_r, t \in \mathcal{T}_i^d \\
        & F^*_{j, k, t} \in \mathbb{R}, \quad \forall (j, k) \in \mathcal{L}, t \in \mathcal{T}^d_i\\
        & \theta^*_{j, t} \in \mathbb{R},  D^*_{j, t} \in \mathbb{R}_+, \quad \forall j \in \mathcal{B}, t \in \mathcal{T}^*_i
        \end{array}\right\}\label{eq:y}
\end{align}
The $i$th DA-UC subproblem is given by \eqref{eq:dauc}
\begin{subequations}\label{eq:dauc}
\begin{align}
    %\begin{split}
        \min &\; f_{u}(\bx^d_i) + f_{e}(\by_i^d) \Delta^d \\
        \textrm{s.t.} &\; c_d(\by^d_i) = 0, \quad c_f(\by^d_i) = 0, \quad c_r(\by^d_i) = 0\\
        &\; c_{su}(\by^d_i, \bx^d_i) \le 0, \quad c_{su, 0}(\by^d_i, \textcolor{blue}{\by^h_k}, \bx^d_i, \bx^d_{i -1}) \le 0\\
        &\;  c_{sd}(\by^d_i, \bx^d_i) \le 0, \quad c_{sd, 0}(\by^d_i, \textcolor{blue}{\by^h_k}, \bx^d_i, \bx^d_{i -1}) \le 0 \\
        &\; x_{g, t} - x_{g, t - 1} = s_{g, t} - z_{g, t} \quad \forall g \in \Gamma_c^d, t \in \mathcal{T}^{d}_i\label{eq:dauc_onoff}\\
        &\; \sum_{j \in {\mathcal{U}^d_{g, t}}} s_{g, j} \le x_{g, t} \quad \forall g \in \Gamma_c^d, t \in \mathcal{T}^d_i \label{eq:dauc_uptime}\\
        &\; \sum_{j \in \mathcal{D}^d_{g, t}}^t z_{g, j} \le 1 - x_{g, t} \quad \forall g \in \Gamma_c^d, t \in \mathcal{T}_i^d \label{eq:dauc_downtime}\\
        &\; \underline{C}_g x_{g, t} \le G^{+, d}_{g, t} + G^{-, d}_{g, t} \le \overline{C}_g x_{g, t}, \quad \substack{\forall g \in \Gamma_c^d, \\t \in \mathcal{T}_i^d} \label{eq:dauc_cap}
    %\end{split}
\end{align}
\end{subequations}
The $i$th subproblem of ST-UC is given by \eqref{eq:stuc}. This has a similar structure as DA-UC; however, there are now DA-UC variables that are incorporated into this lower layer solution through \eqref{eq:stuc_su}, \eqref{eq:stuc_sd}, \eqref{eq:stuc_cap}, and \eqref{eq:stuc_link}. The last constraint ensures that the generation amounts for the DA-UC generators in ST-UC are within a certain bound ($\epsilon^s_g$) of the DA-UC subproblem solutions. This helps avoid myopic solutions, since the ST-UC time horizon is shorter than that of DA-UC. 
\begin{subequations}\label{eq:stuc}
\begin{align}
    %\begin{split}
        \color{red} \min &\; \color{red} f_{u}(\bx^s_i) + f_{e}(\by_i^s) \Delta^s \\
        \color{red} \textrm{s.t.} &\; \color{red} c_d(\by^s_i) = 0, \quad \color{red} c_f(\by^s_i) = 0, \quad \color{red} c_r(\by^s_i) = 0\\
        &\; \color{red} c_{su}(\by^s_i, \textcolor{black}{\bx^d_j}) \le 0, \quad c_{su, 0}(\by^s_i, \textcolor{blue}{\by^h_k}, \textcolor{black}{\bx^d_j}, \textcolor{black}{\bx^d_{j -1}}) \le 0 \label{eq:stuc_su}\\
        &\;  \color{red} c_{sd}(\by^s_i, \textcolor{black}{\bx^d_j}) \le 0, \quad c_{sd, 0}(\by^s_i, \textcolor{blue}{\by^h_k}, \textcolor{black}{\bx^d_j}, \textcolor{black}{\bx^d_{j -1}}) \le 0 \label{eq:stuc_sd} \\
        &\; \color{red} c_{su}(\by^s_i, \bx^s_i) \le 0, \quad c_{su, 0}(\by^s_i, \textcolor{blue}{\by^h_k}, \bx^s_i, \bx^s_{i - 1}) \le 0\\
        &\; \color{red} c_{sd}(\by^s_i, \bx^s_i) \le 0, \quad c_{sd, 0}(\by^s_i, \textcolor{blue}{\by^h_k}, \bx^s_i, \bx^s_{i - 1}) \le 0\\
        &\; \color{red} c_{ru}(\by^s_i) \le 0, \quad \quad \color{red} c_{rd}(\by^s_i) \le 0 \\
        &\; \color{red} x_{g, t} - x_{g, t - 1} = s_{g, t} - z_{g, t} \quad \forall g \in \Gamma_c^s, t \in \mathcal{T}^{s}_i\label{eq:stuc_binary}\\
        &\; \color{red} \sum_{j \in \mathcal{U}^s_{g, t}} s_{g, j} \le x_{g, t} \quad \forall g \in \Gamma_c^s, t \in \mathcal{T}^s_i \label{eq:stuc_ut}\\
        &\; \color{red} \sum_{j \in \mathcal{D}^s_{g, t}} z_{g, j} \le 1 - x_{g, t} \quad \forall g \in \Gamma_c^s, t \in \mathcal{T}_i^s \label{eq:stuc_dt} \\
        &\; \color{red} \underline{C}_g \textcolor{black}{x_{g, \lceil t \rceil}} \le G^{+, s}_{g, t} + G^{-, s}_{g, t} \le \overline{C}_g\textcolor{black}{x_{g, \lceil t \rceil }} \quad \substack{ \forall g \in \textcolor{black}{\Gamma_c^d}, \\ t \in \mathcal{T}_i^s} \label{eq:stuc_cap} \\
        &\; \color{red} \underline{C}_g x_{g, t} \le G^{+, s}_{g, t} + G^{-, s}_{g, t} \le \overline{C}_g x_{g, t} \quad \substack{\forall g \in \Gamma_c^s,\\ t \in \mathcal{T}_i^s}\\
        &\; \color{red} |G^{+, s}_{g, t} + G^{-, s}_{g, t} - \textcolor{black}{G^{+, d}_{g, \lceil t \rceil}} - \textcolor{black}{G^{-, d}_{g, \lceil t \rceil}}| \le \epsilon^s_g, \quad \substack{\forall g \in \textcolor{black}{\Gamma_c^d}, \\t \in \mathcal{T}_i^s} \label{eq:stuc_link}\\
    %\end{split}
\end{align}
\end{subequations}
The $i$th subproblem of HA-ED is given by \eqref{eq:haed}. This formulation is similar to ST-UC but without the binary variables and their accompanying constraints and objective function. In addition, there are now links between {\it both} the ST-UC and DA-UC layers. 
\begin{subequations}\label{eq:haed}
\begin{align}
    %\begin{split}
        \color{blue} \min &\; \color{blue} f_{e}(\by_i^h) \Delta^h \\
        \color{blue} \textrm{s.t.} &\; \color{blue} c_d(\by^h_i) = 0, \quad c_f(\by^h_i) = 0, \quad  c_r(\by^h_i) = 0\\
        &\; \color{blue} c_{su}(\by^h_i, \textcolor{black}{\bx^d_j}) \le 0, \quad c_{su, 0}(\by^h_i, \by^h_{i - 1}, \textcolor{black}{\bx^d_j}, \textcolor{black}{\bx^d_{j - 1}})\\
        &\; \color{blue} c_{sd}(\by^h_i, \textcolor{black}{\bx^d_j}) \le 0, \quad c_{sd, 0}(\by^h_i, \by^h_{i - 1}, \textcolor{black}{\bx^d_j}, \textcolor{black}{\bx^d_{j - 1}})\\
        &\; \color{blue} c_{su}(\by^h_i, \textcolor{red}{\bx^s_i}) \le 0, \quad c_{su, 0}(\by^h_i, \by^h_{i - 1}, \textcolor{red}{\bx^d_j}, \textcolor{red}{\bx^d_{j - 1}})\label{eq:haed_stuc_su} \\
        &\; \color{blue} c_{sd}(\by^h_i, \textcolor{red}{\bx^s_i}) \le 0, \quad c_{sd, 0}(\by^h_i, \by^h_{i - 1}, \textcolor{red}{\bx^d_j}, \textcolor{red}{\bx^d_{j - 1}})\label{eq:haed_stuc_sd} \\
        &\; \color{blue} c_{ru}(\by^h_i) \le 0, \quad \quad c_{rd}(\by^h_i) \le 0 \\
        &\; \color{blue} \underline{C}_g \textcolor{black}{x_{g, \lceil t \rceil}} \le G^{+, h}_{g, t} + G^{-, h}_{g, t} \le \overline{C}_g \textcolor{black}{x_{g, \lceil t \rceil}} \quad \substack{\forall g \in \textcolor{black}{\Gamma_c^d}, \\t \in \mathcal{T}_i^h}\\
        &\; \color{blue} \underline{C}_g \textcolor{red}{x_{g, t}} \le G^{+, h}_{g, t} + G^{-, h}_{g, t} \le \overline{C}_g \textcolor{red}{x_{g, t}} \quad \substack{\forall g \in \textcolor{red}{\Gamma_c^s}, \\t \in \mathcal{T}_i^h}\label{eq:haed_stuc_cap} \\
        &\; \color{blue} |G^{+, h}_{g, t} + G^{-, h}_{g, t} - \textcolor{black}{G^{+, d}_{g, \lceil t \rceil}} - \textcolor{black}{G^{-, d}_{g, \lceil t \rceil}}| \le \epsilon_g, \quad \substack{\forall g \in \textcolor{black}{\Gamma_c^d}, \\t \in \mathcal{T}_i^h } \label{eq:haed_dauclinks}\\
        &\; \color{blue} |G^{+, h}_{g, t} + G^{-, h}_{g, t} - \textcolor{red}{G^{+, s}_{g, \lceil t \rceil}} - \textcolor{red}{G^{-, s}_{g, t}}| \le \epsilon_g, \quad \forall \substack{g \in \textcolor{red}{\Gamma_c^s}, \\t \in \mathcal{T}_i^h} \label{eq:haed_stuc_links}.
    %\end{split}
\end{align}
\end{subequations}

\subsection{Graph Representation}
\begin{figure}[!htp]
    \centering
    \vspace{-0.2in}
    \includegraphics[width=2.0in]{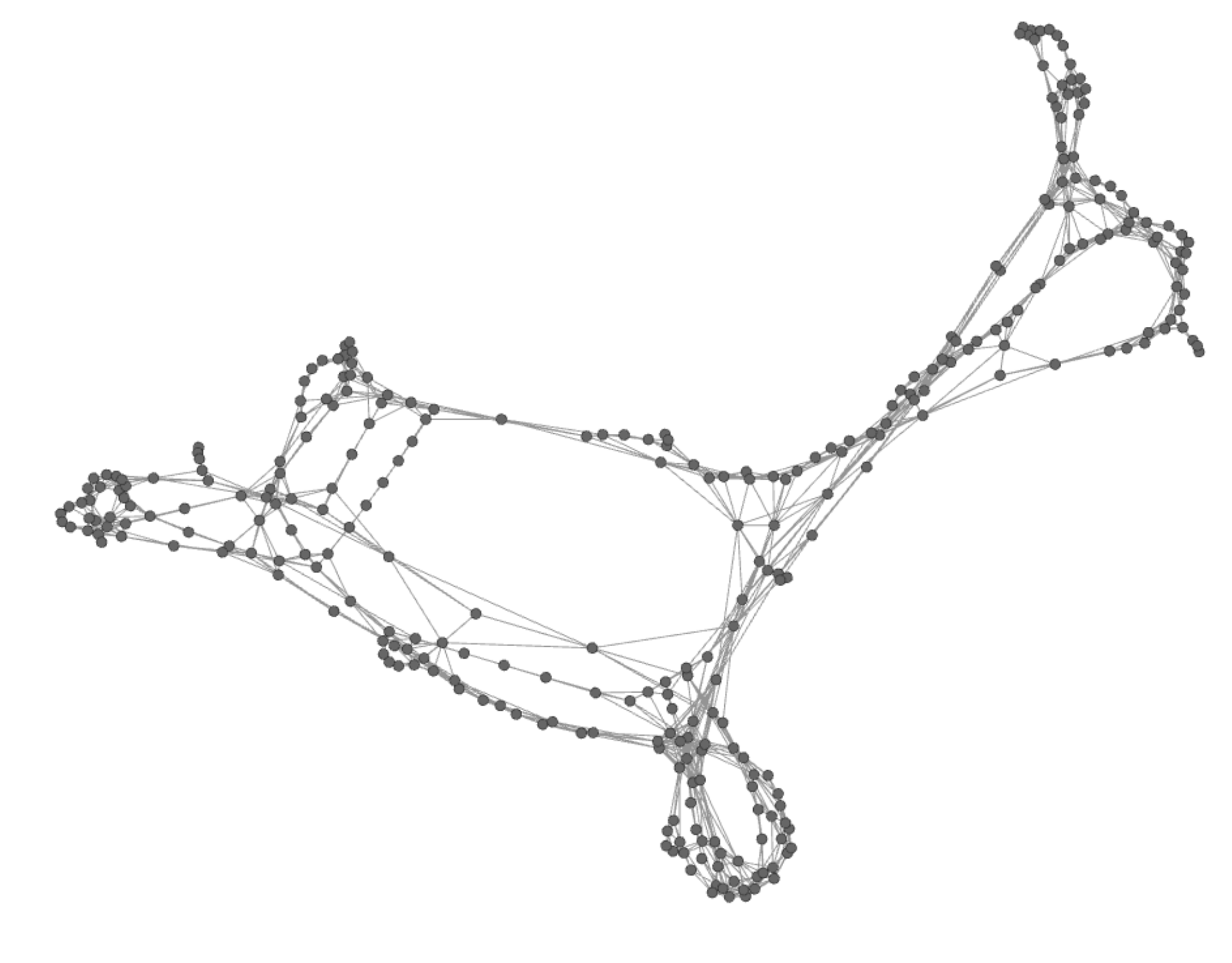}
    \vspace{-0.1in}
    \caption{A graph visualization of the 118-bus system, where nodes correspond to buses and to transmission lines.}
    \label{fig:118_bus_graph}    
\end{figure}

\begin{figure}[!htp]
    \centering
    \vspace{-0.2in}
    \includegraphics[width = 3.25in]{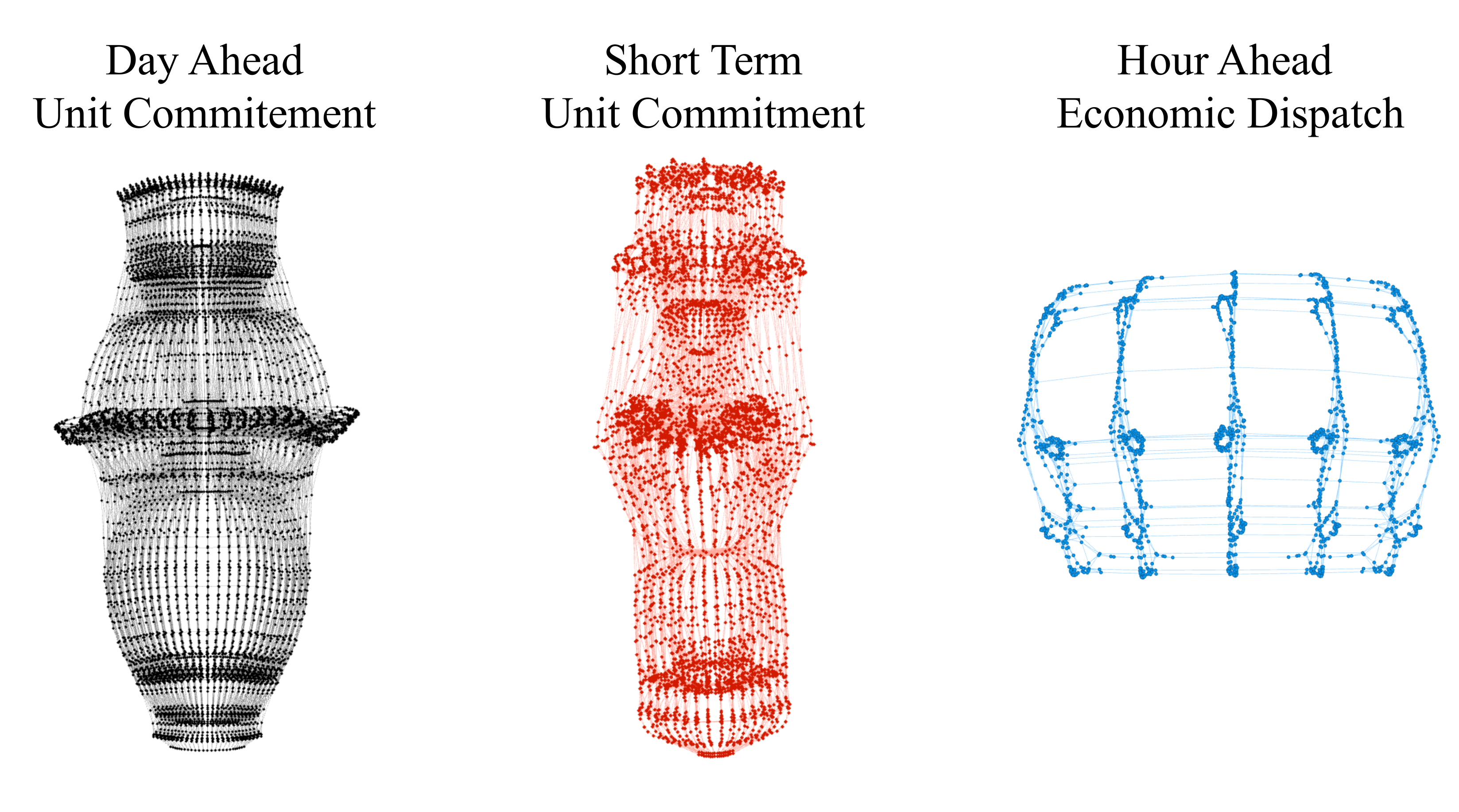}
        \vspace{-0.1in}
    \caption{Representation of a single subproblem subgraph for the DA-UC, ST-UC and HA-ED subproblems.}
    \label{fig:subproblems}
\end{figure}

We now outline how we represent the tri-level hierarchical architecture in {\tt Plasmo.jl}. We represent each time point $t$ as a subgraph, and nodes are placed on this subgraph for each bus and each transmission line. The nodes corresponding to buses contain the variables $D^*_{k, t}$, $\theta^*_{k, t}$, $G^{+, *}_{g, t}$, and $G^{-, *}_{g, t}, \forall g \in \Gamma_k \cap \{ \Gamma_c^* \cup \Gamma_r \}$. In the case of DA-UC and ST-UC problems, the bus nodes also contain $x_{g, t}, s_{g, t}$, and  $z_{g, t} \: \forall g \in \Gamma_k \cap \Gamma_c^*$ and any constraints for these variables. The nodes corresponding to transmission lines contain variables $F^*_{i, j, t}, \forall (i, j) \in \mathcal{L}$. For each node representing line $(i, j) \in \mathcal{L}$, edges (linking constraints) are also placed connecting to bus $i$ and to bus $j$. The resulting subgraph is shown in Figure \ref{fig:118_bus_graph}.  The DA-UC, ST-UC, and HA-ED subproblems were constructed from these time point subgraphs. The DA-UC subgraph has 24 time point subgraphs (i.e., 24 replicates of the network shown in Figure \ref{fig:118_bus_graph}) each representing one hour, the ST-UC subgraph had 16 time point subgraphs with each representing 15 minutes, and the HA-ED subgraph had 5 time point subgraphs with each representing 15 minutes. Linking constraints were also placed between time point subgraphs where applicable, such as for $c_{ru}$, $c_{rd}$, \eqref{eq:dauc_onoff} - \eqref{eq:dauc_downtime}, or \eqref{eq:stuc_binary} - \eqref{eq:stuc_dt}. After subproblem subgraphs were created, the subproblems were combined onto another {\tt OptiGraph} corresponding to one day of operation.  A single day graph contains one DA-UC subgraph, 8 ST-UC subgraphs, and 96 HA-ED subgraphs (105 subgraphs in total). Figure \ref{fig:subproblems} shows the complexity of the resulting graph.

\begin{figure}[!ht]
    \centering
        \vspace{-0.2in}
    \includegraphics[width = 2.0in]{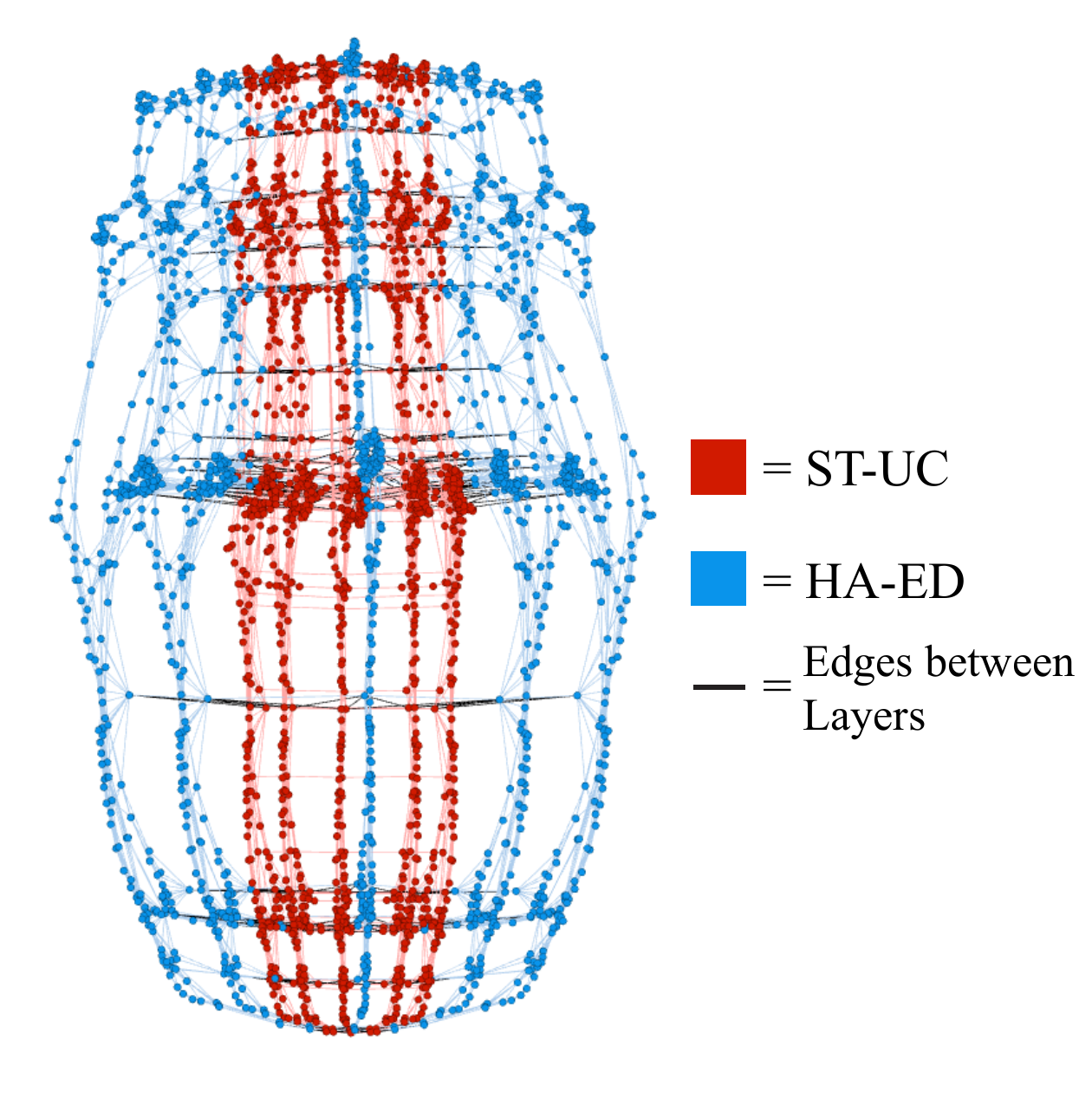}
        \vspace{-0.1in}
    \caption{An example of the hierarchical linking between subproblems. Five time points from an ST-UC subproblem are shown linked to a full HA-ED subproblem, with linking constraints highlighted in black.}
    \label{fig:linking_constraints}
\end{figure}

\begin{figure*}[!t]
    \centering
    \vspace{-0.3in}
    \begin{subfloat}[1-Day Monolith \\ (All Subproblems)]{\includegraphics[width = .3\textwidth]{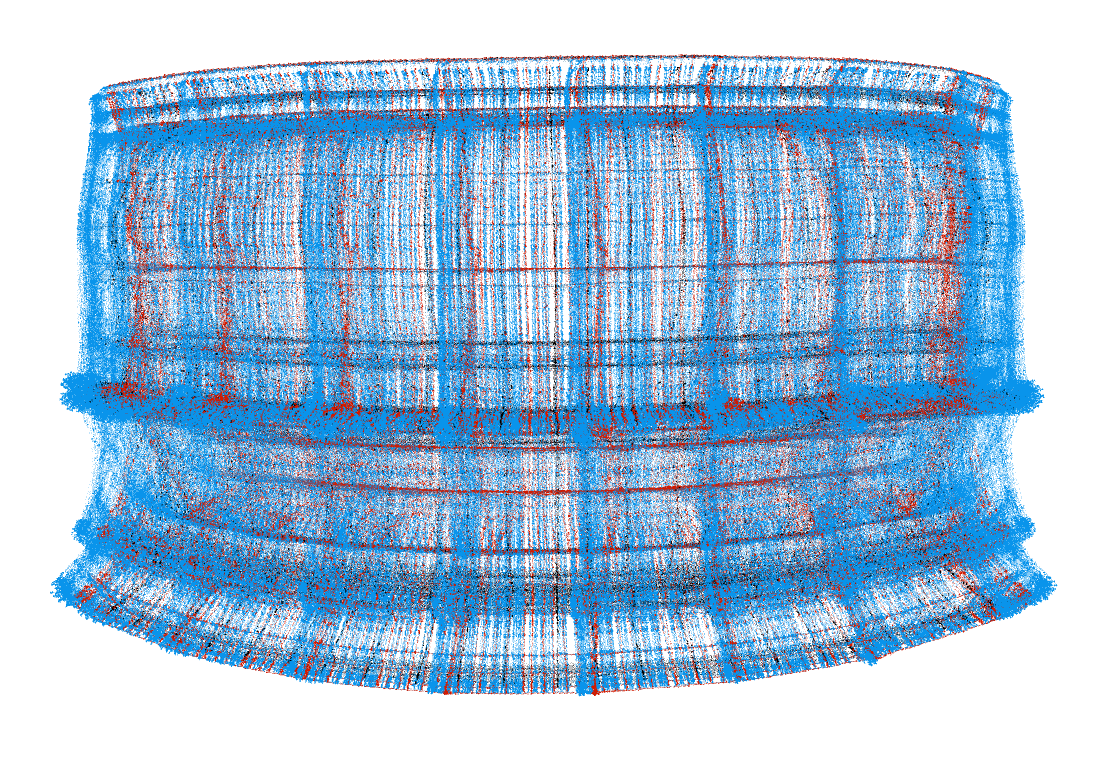}}
    \end{subfloat}
    \hfill
    \begin{subfloat}[Time Point\\ Aggregation]{\includegraphics[width = .25\textwidth]{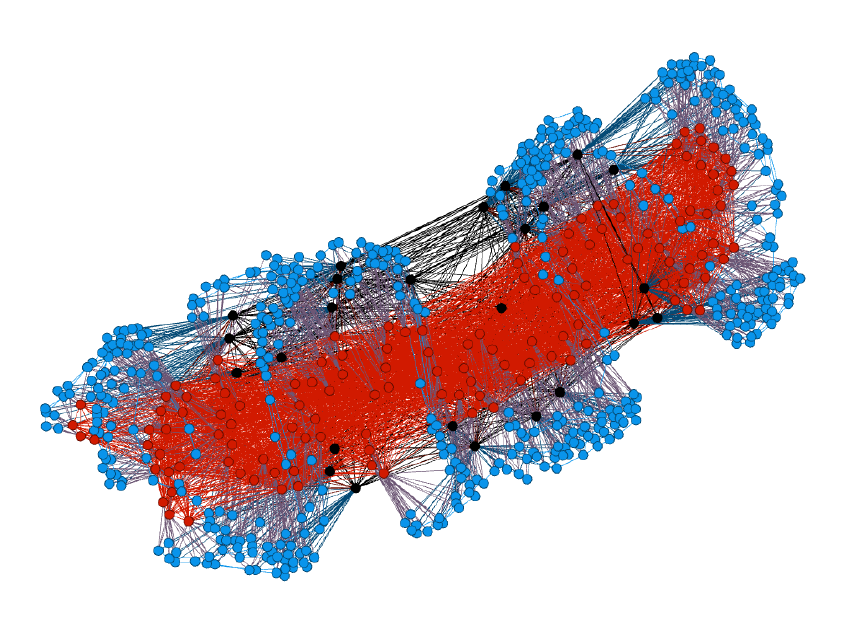}}
    \end{subfloat}
    \hfill
    \begin{subfloat}[Subproblem \\Aggregation]{\includegraphics[width = .2\textwidth]{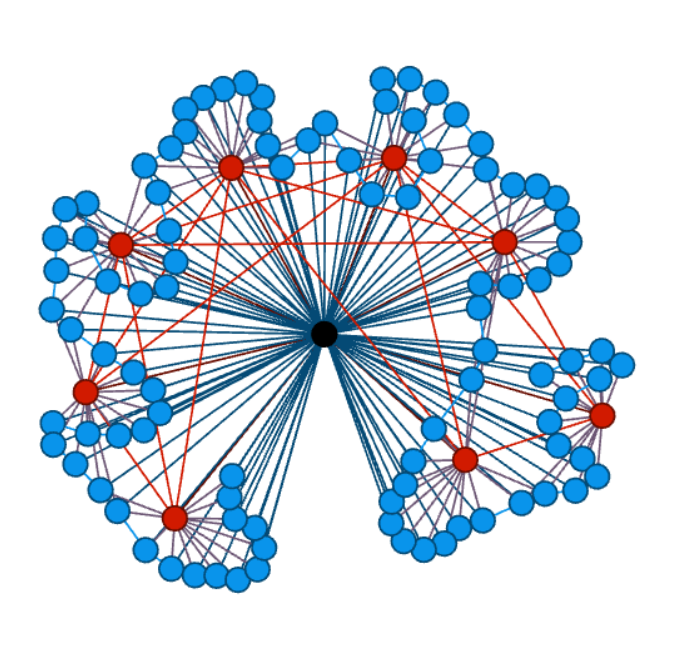}}
    \end{subfloat}
    \begin{subfloat}{\includegraphics[width = .12\textwidth]{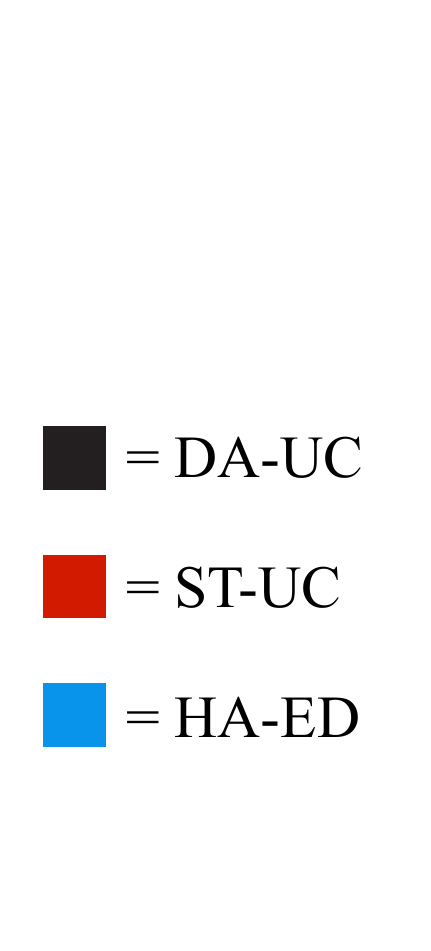}}
    \end{subfloat}
    \caption{Representation of (a) the 1-day monolithic graph highlighted by the subproblems (DA-UC, ST-UC, HA-ED), and equivalent representations where (b) individual time point subgraphs are aggregated into single nodes and where (c) subproblem subgraphs are aggregated into single nodes.}
    \label{fig:aggregation}
            \vspace{-0.2in}
\end{figure*}

With the single-day graph formed, subproblems are linked together according to the formulations given in \eqref{eq:dauc}, \eqref{eq:stuc}, and \eqref{eq:haed}. Figure \ref{fig:linking_constraints} shows an example of this linking for part of an ST-UC subproblem and one HA-ED subproblem. The linking constraints are highlighted in black; these constraints correspond to $c_{su}$ and $c_{sd}$ in \eqref{eq:haed_stuc_su} and \eqref{eq:haed_stuc_sd} and to the linking constraints in \eqref{eq:haed_stuc_cap} and \eqref{eq:haed_stuc_links}. The full problem graph is shown in Figure \ref{fig:aggregation}a, with accompanying representations highlighting the hierarchical structure. The full graph contains 192,432 {\tt OptiNodes} and 292,587 {\tt OptiEdges}. Subgraphs can be collapsed or aggregated into {\tt OptiNodes} without changing the problem formulation and this facilitates visualization. Figure \ref{fig:aggregation}b shows the graph with all time subgraphs aggregated into nodes. Figure \ref{fig:aggregation}c shows all subproblem subgraphs aggregated into nodes; this reveals the hierarchical structure and the linking between layers, where the central black node corresponds to the single-day DA-UC subproblem (top layer), the red nodes correspond to the 8 ST-UC subproblems (middle layer), and the blue nodes correspond to the 96 HA-ED subproblems (bottom layer). 

%\begin{figure*}[!t]
%    \centering
%    \includegraphics[width = 0.9\textwidth]{./figures/subproblem_aggregation.pdf}
%        \vspace{-0.2in}
%    \caption{Representation of the 1-day monolithic graph highlighted by the subproblem (DA-UC, ST-UC, HA-ED). Equivalent representations are shown where individual time point subgraphs are aggregated into single nodes and where subproblem subgraphs are aggregated into single nodes.}
%    \label{fig:aggregation}
%\end{figure*}
%

\subsection{Decomposition Approaches}

Graph representations facilitate the implementation of different decomposition approaches. For example, \cite{atakan2022} decomposed the hierarchy by solving the subproblems in each layer in series (in a receding horizon approach). This sequential decomposition approach can be easily visualized using graphs. For Figure \ref{fig:aggregation}c, this is equivalent to solving the central DA-UC node, then solving the first ST-UC node, and then solving the 12 connected HA-ED nodes. The next ST-UC node is then solved followed by its 12 connected HA-ED nodes (again in order) and so forth until all 8 ST-UC subproblems and their corresponding HA-ED nodes are solved. The solutions of these problems are then passed to the next 1-day graph and the process is repeated. However, there are other decomposition approaches that could be used; for example, instead of solving each subproblem in a receding horizon approach, we could solve each time-point subgraph in a receding horizon approach. This would result in much smaller optimization problems, but likely worse economic performance. In contrast, we could instead solve the entire 1-day monolithic problem as a single optimization problem rather than solving each supbroblem one at a time. The implementation of these strategies can help study trade-offs between tractability and performance. 

\subsection{Results}

\begin{figure}[!ht]
    \centering
            \vspace{-0.0in}
    \includegraphics[width=2.5in]{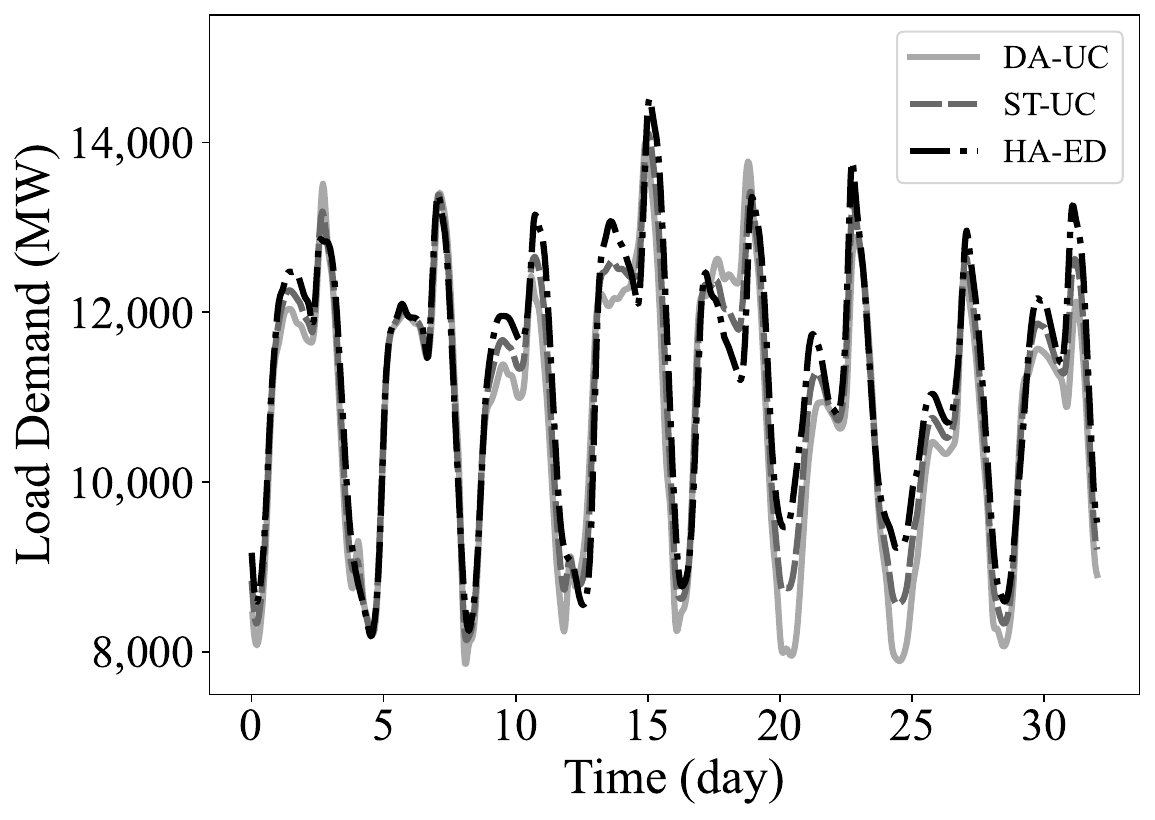}
        \vspace{-0.1in}
    \caption{Load demand used in the DA-UC, ST-UC, and HA-ED subproblems for the 32-day horizon.}
    \label{fig:demands}    
\end{figure}

In this section, we present the results for two different solution approaches. The first approach is to solve the subproblems in a receding horizon approach as done in \cite{atakan2022} ("Receding Horizon"). The second approach is to solve each 1-day graph as a single, monolithic optimization problem ("Monolithic"). Because the monlithic problem is a very large mixed-integer problem (MIP), it took hours to solve, so we used a MIP gap termination criteria of 5\%. In contrast, we used a MIP gap of 0.5\% or less for the receding horizon MIPs as they were smaller and faster to solve. The 1-day monolithic graph contained 641,709 variables (41,727 binary) and 1,103,654 constraints. The code for reproducing these results can be found at \url{https://github.com/zavalab/JuliaBox/tree/master/hierarchical_graphs}. We used the data provided by the 118-bus case study \cite{pena2017extended} used in \cite{atakan2022}. This included a day-ahead (forecasted) load demand, and a real time realized load demand. We used the day-ahead demand for the DA-UC subproblems and the real time demand for the HA-ED subproblems. Because there is no intermediate ``short term" demand data, we used the average of the day ahead and real time demands for the ST-UC subproblems. The three load demands are shown in Figure \ref{fig:demands}. This data was on an hour resolution, so we interpolated the data for higher resolutions. In addition, reserve requirements can vary by system operator, but we chose to use 10\% of the demand for the reserve requirement for UC subproblems and 2.5\% of the demand for ED subproblems which corresponds to the ``low reserve requirements" scenario in \cite{atakan2022}. 

\begin{figure}[!ht]
    \centering
            \vspace{-0.1in}
    \includegraphics[width=3.25in]{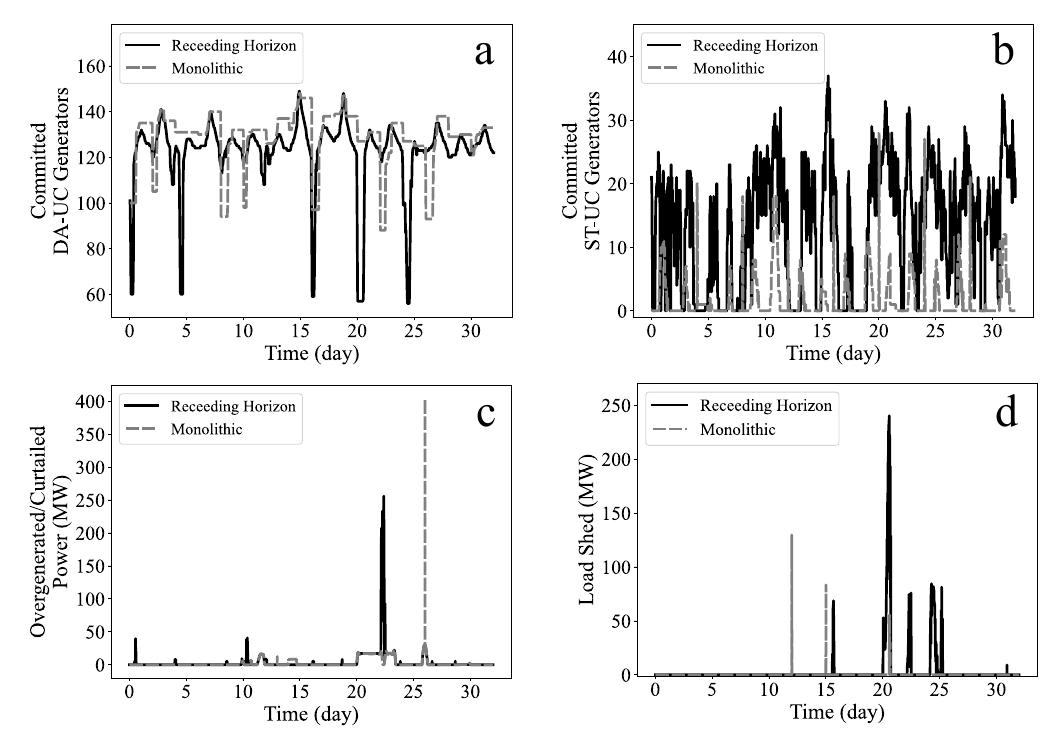}
        \vspace{-0.1in}
    \caption{Results of the tri-level problem using a "receding horizon" type approach and a "monolothic" approach. a) number of committed DA-UC generators over time, b) number of committed ST-UC generators over time, d) overgenerated or curtailed power, and d) load shedding.}
    \label{fig:results}    
\end{figure}

The results of the receding horizon and monolithic approaches are shown in Figure \ref{fig:results} which includes the number of committed DA-UC generators (a), the number of committed ST-UC generators (b), the overgenerated or curtailed power (c), and the amount of load shed (d). The overgenerated/curtailed power and the load shed shown are from the first time point of each HA-ED subproblem. As each HA-ED subproblem had significant overlap with the next problem, we only consider the first HA-ED time point (this is the realized operation). The overall cost of economic dispatch (based on the first time point of each HA-ED subproblem) was \$ 219.4 million  and \$ 199.7 million for the receding horizon and monolithic approaches, respectively.

\subsection{Discussion}

By constructing the hierarchical architecture as a graph, different solution schemes were enabled which provide different insights into the problem. Despite the higher MIP gap used for the monolithic approach, it still performed better than the receding horizon problem and had a lower cost in the economic dispatch by more than \$19 million. The monolithic approach was expected to perform better as the lower layers and upper layers are solved in the same problem, allowing the performance of the lower layer to inform the upper layers. This is also likely why the monolithic approach has less load shedding compared with the receding horizon approach (83.6 MWhr compared with 3113.5 MWhr). The monolithic approach did have a very large peak of overgenerated/curtailed power, but the cost of load shedding (using the costs from \cite{atakan2022}) was 200 times more than the cost of overgenerated/curtailed power. In addition, the monolithic approach had less fluctuation in the number of generators turned on or off.

The results on load shedding were dependent on the reserve requirements used. In this case, the higher reserve requirements on the UC layers compared with the ED layer (10\% vs. 2.5 \%) reduced some of the apparent differences between the demand in the HA-ED layer and the DA-UC layer (e.g., the gap between demand in the DA-UC and the HA-ED layers in day 20 in Figure \ref{fig:demands} would be reduced). If we adjust the reserve requirements and use the ``very low reserve requirements" scenario from  \cite{atakan2022} (5\% of load for UC layer and 1.25\% of load for ED layer), the load shed in the serial problem increases by more than 10 times. While not tested, it is possible that further increasing the reserve requirements could reduce load shedding and/or overgeneration. 

As expected, the monolithic approach took much longer to solve than the receding horizon decomposition approach. In addition, the monolithic approach experienced complications with memory management in the MIP solver. These computational issues, combined with the performance comparisons between the receding horizon and monolithic approaches, highlight the need for decomposition schemes. The receding horizon approach results in a suboptimal solution, but it could be possible to use a decomposition scheme that gives results closer to the monolithic approach but with the computational performance closer to that of the receding horizon problem. Constructing these problems as graphs provides a framework under which a decomposition scheme could be optimized. Overall, this work highlights the utility of representing hierarchical optimization problems using graphs. These graph representations provide a modular way to construct complex (but structured) problems. Each time point can be constructed in a modular manner, and then each time point can be embedded to a modular representation of each subproblem. Graphs are also intuitive to visualize, potentially leading to insights into the problem structure. They provide a framework for manipulating problem structure, such as partitioning/aggregating subgraphs. Graphs also provide a structure that could be exploited via decomposition schemes such as Benders decomposition and Lagrangian relaxation. 

\section{Conclusions and Future Work}
We discussed how hierarchical optimization problems can be represented with graphs. We used the package {\tt Plasmo.jl} to build a tri-level hierarchical optimization problem arising in market operations and presented different approaches to solve the problem. We presented visualizations of these graph representations in {\tt Plasmo.jl}, and we presented the results of the two solution approaches. As part of future work, we are interested in using the graph representation for applying and combining decomposition schemes (e.g., Lagrangian decomposition, Benders decomposition, dual dynamic integer programming) to solve large-scale problem instances.  

\section{Acknowledgements}
This work was supported by the U.S. Department of Energy under grant DE-0002722.

\bibliography{IEEEabrv, bibliography}

\end{document}